\documentclass[a4paper,10pt]{article}
\usepackage[dvips]{graphicx}
\usepackage{amsmath}
\usepackage{a4wide}
\usepackage{amssymb}
\usepackage{amsfonts}
\usepackage{epsfig}
\usepackage{color}
\usepackage{enumitem}

\newcommand{\M}{\mathcal{M}}
\newenvironment{definition}{\refstepcounter{NNEX} \noindent \textbf{Definition \arabic{section}.\arabic{NNEX}:}}{}
\newenvironment{proposition}{\refstepcounter{NNEX} \noindent \textbf{Proposition \arabic{section}.\arabic{NNEX}:}}{}
\newenvironment{proof}{\noindent Proof idea:}{\hfill $\square$\\}
\newcommand{\rsa}{\rightsquigarrow}
\newcommand{\G}{\mathbb{G}}
\newcommand{\Gr}{Grass(p,n)}

\title{Synchronization on the circle}
\author{Alain Sarlette, Rodolphe Sepulchre}

\begin{document}
\newcounter{NNEX}[section]
\definecolor{mygray}{rgb}{0.8,0.8,0.8}
\maketitle

\begin{abstract}
The goal of the present paper is to highlight the fundamental differences of so-called \emph{synchronization} or \emph{consensus} algorithms when the agents to synchronize evolve on a \emph{compact homogeneous manifold} (like the circle, sphere or the group of rotation matrices), instead of a vector space. For the benefit of understanding, the discussion is restricted to the circle. First, a fundamental consensus algorithm on $\mathbb{R}^n$ is reviewed, from which continuous- and discrete-time synchronization algorithms on the circle are deduced by analogy. It is shown how they are connected to Kuramoto and Vicsek models from the literature. Their convergence properties are similar to vector spaces only for specific graphs or if the agents are all located within a semicircle. Examples are proposed to illustrate several other possible behaviors. Finally, new algorithms are proposed to recover (almost-)global synchronization properties.
\end{abstract}

\section{Introduction}\label{s0}

During the last decades, \emph{synchronization} phenomena have attracted the attention of numerous researchers from various fields. Many recent engineering applications consider swarms of agents that combine their efforts in a coordinated (``synchronized'') way to achieve a common task, e.g. distributed exploration \cite{BULLO,LEONARDgen} or interferometry \cite{EADSDarwinControl,Interferometry}. From a modeling and analysis viewpoint, collective phenomena are studied by both the theoretically oriented communities of dynamical systems and statistical physics, as well as more practically oriented experimental physics and biology communities (see e.g. \cite{SmaleEmergence,Kuramoto,Sync,VICSEK}). Celebrated examples of simplified models for the study of coordination phenomena include Kuramoto and Vicsek models.

Kuramoto model is proposed in \cite{Kuramoto} to describe the continuous-time evolution of \emph{phase variables} in a population of weakly coupked oscillators. Each agent $k$ is considered as a periodic oscillator of natural frequency $\omega_k$, whose phase --- i.e. position on its cycle --- at time $t$ is $\theta_k(t) \in S^1$. In addition to its natural evolution, each agent is coupled to all the others:
\begin{equation}\label{eq:3:Kuramotomodel}
\tfrac{d}{dt}\theta_k = \omega_k + \sum_{k=1}^N \, \sin(\theta_j - \theta_k) \; , \quad k=1,2,...,N \, .
\end{equation}

Vicsek model is proposed in \cite{VICSEK} to describe the discrete-time evolution of interacting particles that move with unit velocity in the plane. In the absence of noise, the model writes
\begin{eqnarray}
\label{VicsekLaw1} r_k(t+1) & = & r_k(t) + e^{i \theta_k(t)}\\
\label{VicsekLaw2} \theta_k(t+1) & = & \mathrm{arg}\left( \sum_{j \in n_k(t)} \, e^{i \theta_j(t)} \; + e^{i \theta_k(t)} \right) \; , \quad k=1,2,...,N \, ,
\end{eqnarray}
where $\theta_k$ denotes the heading angle of particle $k$ and $r_k$ its position in the (complex) plane. The set $n_k(t)$ is defined to contain all agents $j$ for which $\Vert r_k(t) - r_j(t) \Vert \leq R$ for some fixed $R > 0$, such that agent $k$ is only influenced by ``close enough'' fellows\footnote{Notation $\Vert z \Vert$ denotes the complex norm of $z \in \mathbb{C}$, such that $z = \Vert z \Vert \,e^{i\, \mathrm{arg}(z)}$. The present paper uses the convention that $\mathrm{arg}(0)$ can take any value on $S^1$. It canonically identifies $\mathbb{C} \cong \mathbb{R}^2$, yielding $\Vert z \Vert = \sqrt{z^T z}$ for $z \in \mathbb{R}^2$.}.

In the engineering community, synchronization has focused on several basic problems. One of them is for a set of agents to reach agreement on some ``quantity'' --- e.g. a position, velocity,\dots --- while exchanging information along a limited set of communication links. This so-called \emph{consensus problem} has been extensively studied for quantities on the \emph{vector space} $\mathbb{R}^n$, see e.g. \cite{hendrickx1,MOREAU2,MOREAU,olfati,Tsitsiklis3,TsitsiklisThesis} and \cite{ConsensusReview} for a review; related results are briefly explained in the next section. However, it appears that many interesting engineering applications and phenomenological models, e.g. as presented above, involve nonlinear \emph{manifolds}; this includes distributed exploration of a planet (sphere $S^2$), rigid body orientations (special orthogonal groups $SO(2)$ or $SO(3)$), and collective motion on limit cycles (abstracted as a circle $S^1$).

Vicsek and Kuramoto models illustrate two general types of problems for synchronization on manifolds. In Vicsek model, stabilizing parallel motion requires equal orientations of the agents, which comes down to reaching \emph{agreement on states} $\theta_k$ on the circle. In Kuramoto model, moving with fixed relative phases corresponds to reaching \emph{agreement on a motion} on the circle. On multi-dimensional manifolds, coordinated motion becomes more complex than on $S^1$ where it simply comes down to reaching agreement on a frequency in $\mathbb{R}$; this is the subject of \cite{MY8}. The present paper focuses on state agreement, more commonly called \emph{synchronization}.\\

\begin{definition}\label{def:1:synchronization}
\emph{A swarm of $N$ agents with states $x_k(t) \in \M$ evolving on a manifold $\M$, $k=1,2,...,N$, is said to asymptotically \emph{synchronize} or \emph{reach synchronization} if}
$$x_1 = x_2 = ... = x_N \quad \text{\emph{asymptotically.}}$$
\emph{The asymptotic value of the $x_k$ is called the \emph{consensus value}.}\\
\end{definition}

The circle is probably the simplest nonlinear manifold to highlight the specifities of synchronization on highly symmetric compact nonlinear spaces. Therefore the present paper focuses on the circle and mentions extensions to manifolds at the end. An important aspect is to maintain the \emph{symmetry} of the synchronization problem: the behavior of the swarm must be invariant with respect to a common translation of all the agents (i.e. on $S^1$, under the transformation $\theta_k \rightarrow \theta_k+a$ $\forall k \in \mathcal{V}$); therefore extensions consider manifolds on which ``all points are equivalent'' --- formally, \emph{compact homogeneous manifolds} like $S^n$, $SO(n)$ or the Grassmann manifolds. In order to concentrate on fundamental issues of synchronization and geometry, system dynamics are simplified to first-order integrators.\\

The paper is organized as follows. Section \ref{s2} briefly reviews a fundamental algorithm for synchronization on vector spaces, also known as a ``consensus algorithm'', in continuous- and discrete-time. Section \ref{s3} extends this algorithm to the circle and highlights its link with Kuramoto \cite{Kuramoto} and Vicsek \cite{VICSEK} models. Section \ref{s4} reviews the positive synchronization properties of these algorithms. Section \ref{s5} illustrates the fact that, in contrast to vector spaces, convergence to non synchronized behavior is possible on the circle. Section \ref{s6} presents three recently proposed algorithms for (almost-)global synchronization on the circle: a \emph{modified coupling function}, a \emph{gossip algorithm} where agents randomly select or discard information from their neighbors, and an algorithm based on auxiliary variables. Section \ref{s7} briefly mentions extensions of the framework to compact homogeneous manifolds. General notations and background information about graphs can be found in the Appendix. Proofs are summarized to their main ideas; complete versions can be found in corresponding references and in \cite{MYthesis}.

\section{Consensus algorithms on vector spaces}\label{s2}

The study of synchronization on vector spaces is a widely covered subject in the systems and control literature of the last decade. The material in the present Section is a summary of basic results developed in this framework by several authors, including \cite{hendrickx1,MOREAU2,MOREAU,olfati,Tsitsiklis3,TsitsiklisThesis}; see \cite{ConsensusReview} for a review and \cite{MOREAU} for some examples of applications.\\

Consider a \emph{swarm} of $N$ agents with states $x_k \in \mathbb{R}^n$, $k \in \mathcal{V} = \lbrace 1,2,...,N \rbrace$, evolving under continuous-time dynamics
\begin{equation}\label{eq:3:VectDynContT}
\tfrac{d}{dt} x_k(t) = u_k(t) \; ,  \quad k=1,2...,N
\end{equation}
or discrete-time dynamics
\begin{equation}\label{eq:3:VectDynDiscT}
x_k(t+1) = x_k(t) + u_k(t) \; ,  \quad k=1,2...,N
\end{equation}
where $u_k$ is a coupling term. The goal is to design $u_k$ such that the agents asymptotically synchronize in the sense of Definition \ref{s0}.\ref{def:1:synchronization}, with the following restrictions.
\begin{itemize}[topsep=0mm, parsep=0mm, itemsep=0mm]
\item[1.] Communication constraint: $u_k$ may only depend on information concerning agent $k$ and the agents $j \rsa k$ sending information to $k$ according to some imposed communication graph $\G$ (see Appendix).
\item[2.] Configuration space symmetry: the behavior of the coupled swarm must be invariant with respect to uniform translation of all the agents: defining $y_k(0) = x_k(0) + a$ $\forall k \in \mathcal{V}$ for any $a \in \mathbb{R}^n$, it must hold $y_k(t) = x_k(t) + a$ $\forall k \in \mathcal{V}$ and $\forall t \geq 0$. Therefore $u_k$ may only depend on the \emph{relative} positions of the agents, i.e. on $(x_j-x_k)$ for $j \rightsquigarrow k$.
\item[3.] Agent equivalence symmetry: all the agents in the swarm must be treated equivalently. This implies that (i) the form of $u_k$ must be the same $\forall k \in \mathcal{V}$ and (ii) all $j \rightsquigarrow k$ must be treated equivalently in $u_k$.
\end{itemize}
This problem is traditionally called the \emph{consensus problem on a vector space}. On manifolds, the term ``synchronization problem'' is preferred because the term ``consensus'' can be given a particular meaning different from synchronization, see \cite{MY4}. On vector spaces, ``consensus'' as defined in \cite{MY4} is equivalent to synchronization, so both terms can be used interchangeably.\\

The consensus problem on vector spaces is solved by the linear coupling
\begin{equation}\label{eq:3:VectLinearController}
u_k(t) = \alpha \; \sum_{j=1}^N a_{jk}(t) (x_j(t) - x_k(t))  \; ,  \quad k=1,2...,N
\end{equation}
where $a_{jk}$ is the weight of link $j \rsa k$ and $\alpha$ is a positive gain. The intuition behind (\ref{eq:3:VectLinearController}) is that each agent moves towards its neighbors, in agreement with the traditional meaning given to a ``consensus'' process. In continuous-time, the closed-loop system (\ref{eq:3:VectDynContT}),(\ref{eq:3:VectLinearController}) implies that agent $k$ is moving towards the position in $\mathbb{R}^n$ corresponding to the (positively weighted) arithmetic mean of its neighbors, $\tfrac{1}{d^{(i)}_k} \sum_{j \rsa k} a_{jk} x_{j}$, where in-degree $d^{(i)}_k = \sum_{j \rsa k} \; a_{jk}$.
In discrete-time, $\alpha$ must satisfy $\alpha \, d^{(i)}_k(t) \leq b$ for some constant $b<1$, $\forall k \in \mathcal{V}$. Then (\ref{eq:3:VectDynDiscT}),(\ref{eq:3:VectLinearController}) means that the future position of agent $k$ is at the (positively weighted) arithmetic mean $\tfrac{1}{\beta_k+d^{(i)}_k} (\sum_{j \rsa k} a_{jk} x_j + \beta_k x_k)$ of its neighbors $j \rsa k$ and itself, with non-vanishing weight $\beta_k$.

Clearly, (\ref{eq:3:VectLinearController}) satisfies the three constraints mentioned above.

The convergence properties of the linear consensus algorithm on a vector space are well characterized. An extension of the following basic result in the presence of time delays can be found in \cite{olfati}; the present paper does not consider time delays.\\

\begin{proposition}\label{prop:3:Moreau}
\emph{(adapted from \cite{MOREAU2,MOREAU,olfati,TsitsiklisThesis}) Consider a set of $N$ agents evolving on $\mathbb{R}^n$ according to  (continuous-time) (\ref{eq:3:VectDynContT}),(\ref{eq:3:VectLinearController}) with $\alpha > 0$ or according to (discrete-time) (\ref{eq:3:VectDynDiscT}),(\ref{eq:3:VectLinearController}) with $\alpha d^{(i)}_k(t) \in [0,b]$ $\forall t \geq 0$ and $\forall k \in \mathcal{V}$, for some constant $b \in (0,1)$. Then the agents globally and exponentially converge to synchronization at some constant value $\bar{x} \in \mathbb{R}^n$ if and only if the communication among agents is characterized by a (piecewise continuous) $\delta$-digraph which is uniformly connected.}

\emph{If in addition, $\G$ is balanced for all times, then the consensus value is the arithmetic mean of the initial values: $\bar{x} = \tfrac{1}{N} \sum_{k=1}^N x_k(0)$.}
\end{proposition}\vspace{2mm}

\begin{proof}
For the first part, see \cite{MOREAU2,MOREAU}, or equivalently \cite{hendrickx1} or \cite{TsitsiklisThesis}.
For the second part, it is easy to see that for a balanced graph, $\tfrac{1}{N} \sum_{k=1}^N x_k(t)$ is conserved over time. The conclusion is then obtained by comparing its value for $t=0$ and for $t$ going to $+\infty$.
\end{proof}

The proof of Proposition \ref{s2}.\ref{prop:3:Moreau} essentially relies on the \emph{convexity} of the update law: the position of each agent $k$ for $t > \tau$ always lies in the convex hull of the $x_j(\tau)$, $j=1,2,...,N$. The permanent contraction of this convex hull, at some nonzero minimal rate because weights are non-vanishing, allows to conclude that the agents end up at a consensus value. An obvious negative consequence of Proposition \ref{s2}.\ref{prop:3:Moreau} for non-varying $\G$ is that synchronization cannot be reached if $\G$ is not root-connected.\\

If interconnections are not only balanced, but also undirected and fixed, then the linear consensus algorithm is a gradient descent algorithm for the disagreement cost function
\begin{equation}\label{eq:3:Vect.costfunction}
V_{\text{vect}}(x) = \tfrac{1}{2} \sum_{k=1}^N \sum_{j=1}^N a_{jk} \Vert x_j - x_k \Vert^2 = \Vert (B \otimes I_n) x \Vert^2 = x^T (L \otimes I_n) x
\end{equation}
where $\Vert z \Vert$ denotes the Euclidean norm $\sqrt{z^T z}$ of $z \in \mathbb{R}^m$, $B$ and $L$ are the incidence and Laplacian matrices of $\G$ respectively (see Appendix), $x \in \mathbb{R}^{N n}$ denotes the vector whose elements $(k-1)n+1$ to $kn$ contain $x_k$, and $\otimes I_n$ is the Kronecker product by the $n \times n$ identity matrix.

\section{Consensus algorithms on the circle}\label{s3}

Consider a swarm of $N$ agents with states on the circle $S^1$. The global topology of the circle is fundamentally different from vector spaces, because if $\theta_k$ denotes an angular position on the circle, then $\theta_k + 2 \pi = \theta_k$, that is, translations on the circle are defined modulo $2\pi$ because they correspond to rotations. This difference in topology, imposing a non-convex configuration space, fundamentally modifies the synchronization problem.

The synchronization problem on $S^1$ is considered under the same agent dynamics as on $\mathbb{R}^n$, i.e. (\ref{eq:3:VectDynContT}) or (\ref{eq:3:VectDynDiscT}) with $x_k$ replaced by $\theta_k$. However, for the design of $u_k$, the different topology induces different implications of the configuration space symmetry. The behavior of the swarm must (i) be invariant with respect to a uniform translation of all $\theta_k$ and (ii) be invariant with respect to the translation of any single $\theta_k$ by a multiple of $2 \pi$ --- i.e., if $\phi_k(0) = \theta_k(0) + 2 a \pi$ for some $k \in \mathcal{V}$ and $a \in \mathbb{Z}$, and $\phi_j(0) = \theta_j(0)$ $\forall j \neq k$, then it must hold $\phi_k(t) = \theta_k(t) + 2 a \pi$ $\forall t \geq 0$ and $\phi_j(t) = \theta_j(t)$ $\forall j \neq k$ and $\forall t \geq 0$. This implies that $u_k$ may only depend on \emph{$2 \pi$-periodic functions of the relative positions} $(\theta_j-\theta_k)$ of the agents $j \rightsquigarrow k$. The simple linear algorithm (\ref{eq:3:VectLinearController}) does not satisfy the periodicity required for configuration space symmetry, and therefore cannot be used on the circle. It can however be used to derive algorithms for synchronization on $S^1$ that are similar to (\ref{eq:3:VectLinearController}) when all agents are within a small arc of the circle. The discrete-time and continuous-time cases are treated consecutively. Because of the symmetry with respect to uniform translations on $S^1$, the swarm's behavior is entirely characterized by examining the evolution of \emph{relative} positions.
\vspace{2mm}

\begin{definition}\label{def:3:configuration}
\emph{A \emph{configuration} is a particular set of relative positions of the agents. Thus a configuration is equivalent to a point $(\bar{\theta}_1,\bar{\theta}_2,...,\bar{\theta}_N)$ $\in S^1 \times S^1 \times ... \times S^1$ and all the points obtained by its uniform rotations $(\bar{\theta}_1+a,\bar{\theta}_2+a,...,\bar{\theta}_N+a)$ for $a \in S^1$.}
\end{definition}

\subsection{Discrete-time}
\label{ss:I:circle.discretetime}

\begin{figure}[htb]
\begin{center}
\setlength{\unitlength}{1mm}
\begin{picture}(120,50)(0,5)
\put(0,0){\includegraphics[width=60mm]{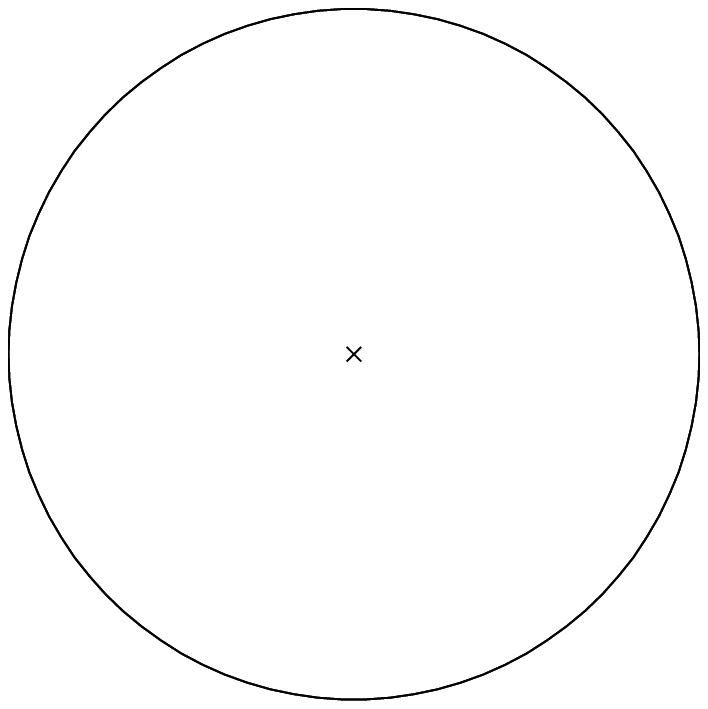}}
\put(0,26.5){\vector(1,0){66}}
\put(60,22){$\Re e$}
\put(31,2){\vector(0,1){50}}
\put(24,50){$\Im m$}
\put(20,44){{\color{blue}\circle*{2}}}
\put(6.5,45){$\theta_{1 \rsa k}(t)$}
\put(20,9){{\color{red}\circle*{2}}}
\put(6.5,6){$\theta_{2 \rsa k}(t)$}
\put(50,36){\circle*{2}}
\put(51.5,33.5){$\theta_k(t)$}
\multiput(31,26.5)(6,3){5}{\line(2,1){4}}
\put(59.5,40.75){\circle{1}}
\multiput(58.5,38.75)(-3,-6){3}{{\color{red}\line(-1,-2){2}}}
\put(50.5,22.75){\circle{1}}
\put(41.5,40.75){\circle{1}}
\multiput(41.5,40.75)(3,-6){3}{{\color{blue}\line(1,-2){2}}}
\put(36,40){$p_k$}
\put(31,26.5){\line(3,4){12}}
\put(43.75,43.5){\circle*{2}}
\put(44,45){$\theta_k(t+1)$}
\put(48.5,37){\vector(-3,4){4}}
\put(47,40){{\footnotesize$u_k(t)$}}
\put(79,12){$p_k := 1.5 e^{i \theta_k} + e^{i \theta_1} + e^{i \theta_2}$}
\end{picture}
\end{center}
\caption{Illustration of update law (\ref{eq:3:circlesynchDT}) for one agent $k$ with $\beta = 1.5$ and $a_{1k} = a_{2k} = 1$.}\label{fig:3:DTcircsynch}
\end{figure}

Synchronization of $\theta_k \in S^1$, $k=1,2,...,N$, can be seen as synchronization of $x_k \in \mathbb{R}^2$ under the constraint $\Vert x_k \Vert = 1$. If the $x_k$ were not restricted to $\Vert x_k \Vert = 1$, algorithm (\ref{eq:3:VectDynDiscT}),(\ref{eq:3:VectLinearController}) would impose $x_k(t+1) = \frac{1}{\beta_k+d^{(i)}_k} \left(\sum_{j = 1}^N a_{jk}\, x_j(t) \; + \beta_k\, x_k(t)\right) \;$, $k=1,2,...,N\;$, with some non-vanishing $\beta_k(t) > 0$. With this update law, $x_k(t+1)$ does generally not satisfy $\Vert x_k(t+1) \Vert = 1$. To obtain $\Vert x_k(t+1) \Vert = 1$, the result of algorithm (\ref{eq:3:VectDynDiscT}),(\ref{eq:3:VectLinearController}) is projected onto the unit circle. Identifying $\mathbb{R}^2 \cong \mathbb{C}$, such that a position on the circle is characterized by $e^{i \theta_k}$, leads to the discrete-time synchronization algorithm
\begin{equation}\label{eq:3:circlesynchDT}
\theta_k(t+1) = \mathrm{arg}\left( \sum_{j = 1}^N a_{jk}\, e^{i \theta_j(t)} \; + \beta \, e^{i \theta_k(t)} \right) \; , \quad k=1,2,...,N \, ,
\end{equation}
for some constant $\beta > 0$. The update of one agent according to (\ref{eq:3:circlesynchDT}) is illustrated on Figure \ref{fig:3:DTcircsynch}. It is clear from the picture that (\ref{eq:3:circlesynchDT}) respects the geometric invariance of $S^1$. This is confirmed by rewriting (\ref{eq:3:circlesynchDT}) as
\begin{equation}\label{eq:3:circlesynchDTalt}
\theta_k(t+1) = \theta_k(t) + u_k = \theta_k(t) + \mathrm{arg}\left( \sum_{j = 1}^N a_{jk}\, e^{i (\theta_j(t)-\theta_k(t))} \; + \beta \right) \; , \quad k=1,2,...,N
\end{equation}
where $u_k$ indeed only involves $2\pi$-periodic functions of relative positions of connected agents $j \rsa k$.

For fixed undirected $\G$, the point $\theta_k(t+1)$ obtained from (\ref{eq:3:circlesynchDT}) is the projection on the unit circle of a point obtained by gradient descent for $V_{\text{vect}}$ in the complex plane.

\paragraph*{Vicsek model:} Heading update law (\ref{VicsekLaw2}) actually corresponds to (\ref{eq:3:circlesynchDT}) with $\beta = 1$ and $a_{jk} \in \lbrace 0,1\rbrace$. The positions define the interconnection graph $\G(t)$ by imposing $a_{jk} = 1$ if and only if $\Vert r_k(t) - r_j(t) \Vert \leq R$; therefore $\G$ is called a \emph{proximity graph}. The study of proximity graphs, or other state-dependent graphs, is beyond the scope of the present work.

\subsection{Continuous-time}
\label{ss:II:circle.continuoustime}

Taking the continuous-time limit of (\ref{eq:3:circlesynchDT}) amounts to letting $\beta$ grow indefinitely. In this limit case, $x_k(t+1) \in \mathbb{C}$ is defined with an infinitesimal gradient step for $V_{\text{vect}}$, and projected onto $S^1$ to yield $\theta_k(t+1)$. This is strictly equivalent to projecting the gradient of $V_{\text{vect}}$ onto the tangent to the circle at $\theta_k(t)$, and taking a corresponding infinitesimal descent step along the circle. Thus by viewing $V_{\text{vect}}$ as a function of $\theta$, renamed for clarity
\begin{equation}\label{eq:3:CostFunctionTheta}
V_{\text{circ}}(\theta) = \tfrac{1}{2} \sum_{k=1}^N \sum_{j=1}^N a_{jk} \Vert e^{i \theta_j} - e^{i \theta_k} \Vert^2 = \tfrac{1}{2} \sum_{k=1}^N \sum_{j=1}^N a_{jk} \left( 2 \sin(\tfrac{\theta_j-\theta_k}{2}) \right)^2\; ,
\end{equation}
the corresponding gradient descent algorithm along the circle, $\tfrac{d}{dt}\theta_k = - \alpha \frac{\partial V_{\text{circ}}}{\partial \theta_k}$ $\forall k \in \mathcal{V}$ with $\alpha > 0$, is the continuous-time limit of (\ref{eq:3:circlesynchDT}). Computing the gradient of (\ref{eq:3:CostFunctionTheta}) yields the following continuous-time algorithm for synchronization on the circle, with constant $\alpha > 0$:
\begin{equation}\label{eq:3:CircleSynchCTbef}
\tfrac{d}{dt} \theta_k = \alpha \; \sum_{j=1}^N (a_{jk}(t)+a_{kj}(t))\; \sin(\theta_j(t) - \theta_k(t)) \; , \quad k=1,2,...,N \; .
\end{equation}
This algorithm can only be implemented for undirected $\G$. An extension to directed graphs is
\begin{equation}\label{eq:3:CircleSynchCT}
\tfrac{d}{dt} \theta_k = 2 \alpha \; \sum_{j=1}^N a_{jk}(t)\; \sin(\theta_j(t) - \theta_k(t)) \; , \quad k=1,2,...,N \, .
\end{equation}
This is the actual algorithm considered in the following. It satisfies the geometric invariance of $S^1$, since the right-hand side is a $2 \pi$-periodic function of relative positions. With $x_k = e^{i \theta_k}$, (\ref{eq:3:CircleSynchCT}) is equivalent to
\begin{equation}\label{eq:3:CircleSynchCTalt}
\tfrac{d}{dt} x_k = 2 \alpha \, \mathrm{Proj}_{x_k} \left( \sum_{j=1}^N a_{jk} \, (x_j - x_k) \right)
\end{equation}
where $\mathrm{Proj}_{x_k}(r_k) = r_k - x_k \, x_k^T r_k$ denotes the orthogonal projection of $r_k \in \mathbb{C} \cong \mathbb{R}^2$ onto the direction tangent to the unit circle at $x_k = e^{i \theta_k}$. The geometric interpretation is that (\ref{eq:3:CircleSynchCTalt}) defines a consensus update similar to (\ref{eq:3:VectDynContT}),(\ref{eq:3:VectLinearController}) but constrained to the subset of $\mathbb{R}^2$ where $\Vert x_k \Vert = 1$. Algorithm (\ref{eq:3:CircleSynchCT}) was proposed in \cite{BlockIsland} in a control framework, and is directly linked to Kuramoto model.

\paragraph*{Kuramoto model:} Comparing with (\ref{eq:3:Kuramotomodel}), algorithm (\ref{eq:3:CircleSynchCT}) in fact corresponds to Kuramoto model with equal natural frequencies $\omega_1=\omega_2=...=\omega_N$, but general interconnections. This highlights a link between the sine-model of Kuramoto and the ``averaging'' update law for headings in Vicsek model. For the complete graph, $V_{\text{circ}} = \tfrac{1}{2} \sum_{k=1}^N \sum_{j=1}^N \Vert e^{i \theta_j} - e^{i \theta_k} \Vert^2 = N^2 - \left\Vert \sum_{k=1}^N e^{i \theta_k} \right\Vert^2 \;$.
The quantity $\left\Vert \sum_{k=1}^N e^{i \theta_k} \right\Vert^2$, known as the ``complex order parameter'' in the context of Kuramoto model, has been used for decades as a measure of the synchrony of phase variables in the literature on coupled oscillators.

The main point in studies of Kuramoto model is the coordination of agents having different $\omega_k$. The important issue of robustly coordinating agents despite their different natural tendencies is not the subject of the present paper.

\section{Convergence properties}\label{s4}

Section \ref{s3} proposes algorithms (\ref{eq:3:circlesynchDT}) and (\ref{eq:3:CircleSynchCT}) as natural extensions of synchronization algorithms for the circle. However, the circle is not a convex configuration space. As a consequence, the convergence properties of (\ref{eq:3:circlesynchDT}) and (\ref{eq:3:CircleSynchCT}) do not match those of (\ref{eq:3:VectDynContT}),(\ref{eq:3:VectLinearController}) and (\ref{eq:3:VectDynDiscT}),(\ref{eq:3:VectLinearController}) on vector spaces. The present Section focuses on positive convergence results, while Section \ref{s5} focuses on situations in which asymptotic synchronization is not achieved.

\subsection{Local synchronization like for vector spaces}

When all agents are within a small subset of $S^1$, (\ref{eq:3:CircleSynchCT}) becomes similar to (\ref{eq:3:VectDynContT}),(\ref{eq:3:VectLinearController}) because $\sin(\theta_j - \theta_k) \simeq (\theta_j - \theta_k)$ for small $(\theta_j - \theta_k)$. A similar observation can be made for the discrete-time algorithms. It is thus not surprising that \cite{jadbabaie3,MOREAU} are able to show that asymptotic synchronization on the circle is locally achieved under the same conditions as on vector spaces.\\

\begin{proposition}\label{prop:3:local}
\emph{(adapted from \cite{MOREAU}) Consider a set of $N$ agents evolving on $S^1$ according to (continuous-time) (\ref{eq:3:CircleSynchCT}) with $\alpha > 0$ or according to (discrete-time) (\ref{eq:3:circlesynchDT}) with $\beta > 0$. If the communication among agents is characterized by a (piecewise continuous) $\delta$-digraph $\G$ which is uniformly connected and all agents are initially located within an open semicircle, then they exponentially converge to synchronization at some constant value $\bar{\theta} \in S^1$.}
\end{proposition}\vspace{2mm}

\begin{proof} 
Assume without loss of generality that $\theta_k \in [-b,b] \subset (-\tfrac{\pi}{2},\tfrac{\pi}{2})$ initially.  Then $\sin(\theta_k(t)-\theta_j(t)) = c_{jk}(t) (\theta_j(t)-\theta_k(t))$ where $c_{jk}(t) \geq \tfrac{\sin(2b)}{2b} > 0$ depends on $(\theta_j(t)-\theta_k(t))$. Thus (\ref{eq:3:CircleSynchCT}) is equivalent to 
\begin{equation}\label{relect1:1:eq}
\tfrac{d}{dt}\theta_k = \alpha \sum_{j=1}^N a_{jk}(t) \, c_{jk}(t)\; (\theta_j(t) - \theta_k(t))
\end{equation}
for some time-varying $c_{jk} \geq \tfrac{\sin(2b)}{2b} > 0$. This can be viewed as a linear synchronization algorithm for $\theta \in \mathbb{R}$ with a $\delta_2$-digraph $\G_2$ of weights $(a_{jk} c_{jk})$. A similar idea can be used in discrete-time.
\end{proof}

When the agents are distributed over more than a semicircle, the proof of Proposition \ref{s4}.\ref{prop:3:local} no longer holds, because the $c_{jk}$ can be negative. This is in fact the consequence of a \emph{loss of convexity}, implying that the strong vector space arguments of \cite{MOREAU2,MOREAU} are no longer applicable\footnote{An open subset $s \subset S^1$ is \emph{convex} if it contains all shortest paths between any two points of $s$.}. In algorithms (\ref{eq:3:circlesynchDT}) and (\ref{eq:3:CircleSynchCT}), agents move on the shortest path towards their neighbors. Therefore, if the agents are initially located within a semicircle, then they remain within this set for all future times, while agents distributed over more than a semicircle may, a priori, leave any open arc $s \in S^1$ containing them all. A more global analysis requires stronger assumptions.

\subsection{Some graphs ensure (almost) global synchronization}

For general graphs, synchronization is only locally asymptotically stable. For some graphs however, synchronization is (almost) globally asymptotically stable.\\

\begin{proposition}\label{prop:3:SynchGraphs}
\emph{(adapted from \cite{KurSynchGraph,SPL2005}) Consider a set of $N$ agents evolving on $S^1$ by applying (\ref{eq:3:CircleSynchCT}) with $\alpha > 0$ or (\ref{eq:3:circlesynchDT}) with $\beta > 0$ not too small (see Propositions \ref{s5}.\ref{prop:3:CircleAsynch} and \ref{s5}.\ref{prop:3:DTSynch} in the following section). If communication graph $\G$ is either a fixed directed root-connected tree,
or an undirected tree, or the complete graph, or any vertex-interconnection of trees and complete graphs\footnote{A vertex-interconnection of two graphs $\G_1(\mathcal{V}_1,\mathcal{E}_1)$ and $\G_2(\mathcal{V}_2,\mathcal{E}_2)$ is a graph $\G$ whose vertices can be partitioned into a singleton $\lbrace k \rbrace$ and two sets $\mathcal{V}_a$, $\mathcal{V}_b$ and whose edge set can be partitioned into two sets $\mathcal{E}_a$, $\mathcal{E}_b$, such that $\mathcal{V}_a \cup \lbrace k \rbrace = \mathcal{V}_1$, $\mathcal{E}_a = \mathcal{E}_1$, $\mathcal{V}_b \cup \lbrace k \rbrace = \mathcal{V}_2$ and $\mathcal{E}_b = \mathcal{E}_2$.}, then the agents asymptotically converge to synchronization, for almost all initial conditions.}
\end{proposition}\vspace{2mm}

\begin{proof}
For the directed rooted tree, each agent is attracted towards its parent and, except for unstable situations where two connected agents are exactly at opposite positions on the circle, they synchronize at the initial position of the root.
The particular undirected graphs have the property that synchronization is the only local minimum of $V_{\mathrm{circ}}$.  This is rather obvious for the tree, from the fact that for any pair of connected agents, variations of $V_{\mathrm{circ}}$ can be built involving only the distance between that particular pair of agents. The property is proved for the complete graph in \cite{SPL2005}. Finally, \cite{KurSynchGraph} shows that the property holds for vertex-interconnections of graphs for which it holds individually. Then Propositions \ref{s5}.\ref{prop:3:cont.undir}, \ref{s5}.\ref{prop:3:CircleAsynch} and \ref{s5}.\ref{prop:3:DTSynch} of the following section ensure that synchronization is the only stable limit set.
\end{proof}

\section{Obstacles to global synchronization}\label{s5}

Section \ref{s4} identifies situations where (\ref{eq:3:circlesynchDT}) and (\ref{eq:3:CircleSynchCT}) converge to synchronization. The present section examines what can happen when this is not the case.

\subsection{Convergence to local equilibria for fixed undirected $\G$}\label{ss:ConvToEq}

The fact that for fixed undirected $\G$, (\ref{eq:3:CircleSynchCT}) is a gradient descent for $V_{\text{circ}}$, has strong implications for the convergence analysis.\vspace{2mm}

\begin{proposition}\label{prop:3:cont.undir}
\emph{Consider a set of $N$ agents evolving on $S^1$ according to (\ref{eq:3:CircleSynchCT}) with $\alpha > 0$, with communication graph $\G$ fixed and undirected. Then the agents always converge to a set of equilibria corresponding to the critical points\footnote{A \emph{critical point} of a differentiable function $f: \mathcal{X} \rightarrow \mathbb{R}$ is a point of $\mathcal{X}$ where the gradient of $f$ is identically zero.} of $V_{\mathrm{circ}}$ defined in (\ref{eq:3:CostFunctionTheta}). The only asymptotically stable equilibria are the local minima of $V_{\mathrm{circ}}$.}
\end{proposition}\vspace{2mm}

\begin{proof}
Properties of gradient algorithms.
\end{proof}

For the discrete-time algorithm, in addition to a result similar to Proposition \ref{s5}.\ref{prop:3:cont.undir} but with a bound on $\beta$, it can be shown that under \emph{locally asynchronous update}, convergence holds for arbitrary $\beta$, i.e. without requiring any minimal inertia. Agents are said to update \emph{synchronously} if, between instants $t$ and $t+1$, all agents $k \in \mathcal{V}$ apply (\ref{eq:3:circlesynchDT}). In contrast, agents are said to update \emph{locally asynchronously} if, between instants $t$ and $t+1$, only a subset of agents $\sigma \subset \mathcal{V}$ applies (\ref{eq:3:circlesynchDT}) and the others remain at their position, and set $\sigma$ contains no agents that are connected to each other in $\G$.\\

\begin{proposition}\label{prop:3:CircleAsynch}
\emph{(from \cite{MY1}) Consider a set of $N$ agents evolving on $S^1$ by applying (\ref{eq:3:circlesynchDT}) locally asynchronously with update subset sequence $\sigma(t)$, for $\beta > 0$ and fixed undirected $\G$. Assume that there exist a finite time span $T$ and a partition of the discrete-time space $[t_0, t_1)\, , \, [t_1 , t_2)$,... with $(t_{n+1} - t_n) < T$ $\forall n \in \mathbb{Z}_{\geq 0}$, such that $\forall k \in \mathcal{V}$ and $\forall n \in \mathbb{Z}_{\geq 0}$, there exists $t \in [t_n, t_{n+1})$ such that $k \in \sigma(t)$. Then the agents almost always converge to a set of equilibria corresponding to the critical points of $V_{\mathrm{circ}}$. The only asymptotically stable equilibria are the local minima of $V_{\mathrm{circ}}$.}
\end{proposition}\vspace{2mm}

\begin{proof}
Under the assumptions of the Proposition, denoting $\sum_{j=1}^N \, a_{jk} e^{i (\theta_j-\theta_k)} + \beta = \rho_k e^{i u_k}$,
$$V_{\mathrm{circ}}(t+1)-V_{\mathrm{circ}}(t) = -2 \, \sum_{k \in \sigma(t)} \; (\rho_k + \beta) \, \left(\sin(\tfrac{u_k}{2})\right)^2 \leq 0 \, .$$
Since every agent $k$ is updated at an infinite number of time instants, $u_k$ or $\rho_k$ (making $u_k$ undefined) must go to $0$ when $t$ goes to $+\infty$, $\forall k \in \mathcal{V}$. The case $\rho_k = 0$ has zero measure and can only appear ``by chance'', because it is the global maximum of $V_{\mathrm{circ}}$ with respect to $\theta_k$. If $u_k$ goes to $0$, then agent $k$ asymptotically approaches an equilibrium set; only minima can be asymptotically stable for a descent algorithm.
\end{proof}

For synchronous update, it is again necessary to impose a bound on the motion of the agents. However this bound is not easy to find. The following result provides a conservative bound on $\beta$.\\

\begin{proposition}\label{prop:3:DTSynch}
\emph{(from \cite{MY1}) Consider a set of $N$ agents evolving on $S^1$ by applying (\ref{eq:3:circlesynchDT}) synchronously with fixed, undirected and unweighted $\G$. Assume that
$$\beta \geq d_{max} \, (\tfrac{2}{M^{\ast}} + 1) \quad \text{ where } \quad \tfrac{e^{M^{\ast}}-1}{M^{\ast}} = 1 + \tfrac{d_{max}}{d_{sum}}$$
with $d_{sum} = \sum_{j=k}^N \, d_k^{(i)}$ and $d_{max} = \mathrm{max}_{k \in \mathcal{V}}(d_k^{(i)})$, where $d_k^{(i)}$ is the in-degree of agent $k$. Then the agents almost always converge to a set of equilibria corresponding to the critical points of $V_{\mathrm{circ}}$. The only asymptotically stable equilibria are the local minima of $V_{\mathrm{circ}}$.}
\end{proposition}\vspace{2mm}

\begin{proof}
The proof shows that $V_{\mathrm{circ}}(t+1)-V_{\mathrm{circ}}(t) \leq 0$ for synchronous operation and the bound on $\beta$. See \cite{MY1} or \cite{MYthesis} for complete computations.
\end{proof}

A problem with the bound of Proposition \ref{s5}.\ref{prop:3:DTSynch} is that each agent must know $d_{sum}$ and $d_{max}$, which is information about the (communication structure of the) whole swarm.

In the absence of inertia ($\beta = 0$), (\ref{eq:3:circlesynchDT}) can lead to a limit cycle in synchronous operation, at least for some $\G$ (see \cite{MYthesis} for an example).

\paragraph*{Hopfield network:} This model proposed in \cite{HOPFIELD} considers $N$ neurons with states $x_k$ $\in$ $\lbrace -1, 1\rbrace$. The discrete-time update law for the states of the neurons is
\begin{equation}\label{eq:3:Hopf.1}
x_k(t+1) = \mathrm{sign}\left( \sum_{j=1}^N \, a_{jk} x_j(t) + \xi_k \right) \; , \quad 1,2,...,N
\end{equation}
where $\xi_k$ is a firing threshold. Considering
$V_H = \frac{-1}{2} \sum_{k=1}^N \sum_{j=1}^N \; a_{jk} \, x_j \, x_k - \sum_{k=1}^N x_k \, \xi_k \; ,$
\cite{HOPFIELD} shows that when (\ref{eq:3:Hopf.1}) is applied asynchronously with a random update sequence,
the property $V_H(t+1) \leq V_H(t)$ always holds and the network eventually reaches a local minimum of $V_H$. In contrast,  the system can go into a limit cycle under synchronous operation (see \cite{GOLES}).

Defining the sphere $S^n$ of dimension $n$ as $\lbrace x_k \in \mathbb{R}^{n+1} : \Vert x_k \Vert = 1 \rbrace$, the set $\lbrace -1, 1 \rbrace$ can be seen as ``$\;S^0\;$'', while the circle is $S^1$. For $\xi_k = 0$, (\ref{eq:3:Hopf.1}) is in fact the strict analog of (\ref{eq:3:circlesynchDT}) for ``the sphere of dimension 0'' --- namely moving towards the neighbors in the embedding vector space and projecting back to the state space. The absence of inertia in (\ref{eq:3:Hopf.1}) would correspond to $\beta=0$ in (\ref{eq:3:circlesynchDT}). Both (\ref{eq:3:circlesynchDT}) and (\ref{eq:3:Hopf.1}) can be viewed as projections of descent algorithms for a symmetric quadratic potential, which remain descent algorithms under locally asynchronous update such that convergence is ensured. Both algorithms can fail to converge and run into a limit cycle in synchronous operation.\\

Propositions \ref{s5}.\ref{prop:3:cont.undir}, \ref{s5}.\ref{prop:3:CircleAsynch} and \ref{s5}.\ref{prop:3:DTSynch} say that the stable equilibria are the minima of $V_{\mathrm{circ}}$. Unfortunately, depending on $\G$, there may be local minima different from synchronization.

\begin{figure}[htb]
\begin{center}\setlength{\unitlength}{1mm}
\begin{picture}(150,50)
\multiput(0,0)(50,0){3}{\includegraphics[width=50mm]{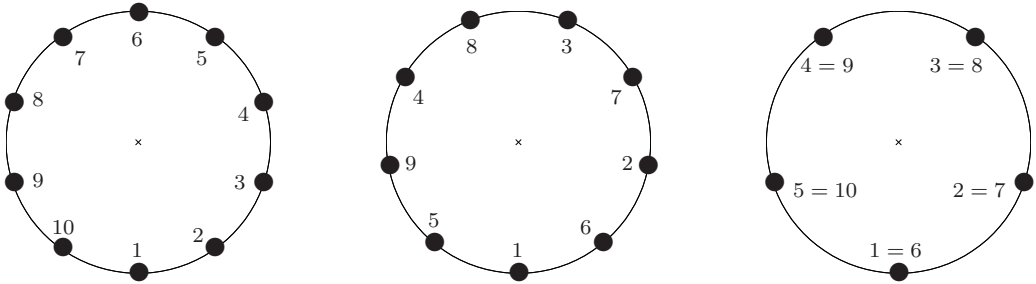}}
\multiput(26,5)(50,0){3}{\circle*{2.5}}
\put(26,39.5){\circle*{2.5}}
\put(36,8.3){\circle*{2.5}}
\put(16,8.3){\circle*{2.5}}
\multiput(36,36.2)(100,0){2}{\circle*{2.5}}
\multiput(16,36.2)(100,0){2}{\circle*{2.5}}
\put(42.4,27.6){\circle*{2.5}}
\put(9.6,27.6){\circle*{2.5}}
\multiput(42.4,16.9)(100,0){2}{\circle*{2.5}}
\multiput(9.6,16.9)(100,0){2}{\circle*{2.5}}
\put(122,7){{\footnotesize $1 = 6$}} \put(133,15){{\footnotesize $2 = 7$}}
\put(130,31.5){{\footnotesize $3 = 8$}} \put(113,31.5){{\footnotesize $4 = 9$}}
\put(112,15){{\footnotesize $5 = 10$}}
\multiput(25,7)(50,0){2}{{\footnotesize $1$}}\put(33,9){{\footnotesize $2$}}\put(38.5,16){{\footnotesize $3$}}
\put(39,25){{\footnotesize $4$}}\put(33.5,32.5){{\footnotesize $5$}}\put(25,35){{\footnotesize $6$}}
\put(17.5,32.5){{\footnotesize $7$}}\put(12,27){{\footnotesize $8$}}\put(12,16.3){{\footnotesize $9$}}
\put(14.5,10){{\footnotesize $10$}}
\multiput(64.9,9)(22.2,0){2}{\circle*{2.5}} \put(64,11){{\footnotesize $5$}}\put(84,10){{\footnotesize $6$}}
\multiput(59,19.2)(34,0){2}{\circle*{2.5}} \put(61,18.5){{\footnotesize $9$}}\put(89.5,18.5){{\footnotesize $2$}}
\multiput(61.05,30.9)(29.9,0){2}{\circle*{2.5}} \put(62,27.3){{\footnotesize $4$}}\put(88,27.3){{\footnotesize $7$}}
\multiput(69.6,38.45)(12.8,0){2}{\circle*{2.5}} \put(69,34){{\footnotesize $8$}}\put(81.5,34){{\footnotesize $3$}}
\end{picture}
\end{center}
\caption{Several balanced configurations that are stable for the undirected ring graph; agents are numbered in the order of the ring, e.g. agent $3$ is connected to agents $2$ and $4$.}\label{fig:Splayring}
\end{figure}

\paragraph*{Local equilibria for the undirected ring:} The following example is taken from \cite{Polygons}. Consider $N$ agents interconnected according to an undirected, unweighted ring graph. Then the critical points of $V_{\mathrm{circ}}$ satisfy $\; \sin(\theta_{j_a(k)}-\theta_k) + \sin(\theta_{j_b(k)}-\theta_k) = 0 \;$ where $j_a(k)$ and $j_b(k)$ are the two neighbors of $k$ in the ring graph. This requires positions of consecutive agents in the ring graph to differ either by $\theta_0$ or by $\pi-\theta_0$, for some $\theta_0 \in [0,\tfrac{\pi}{2}]$ well chosen such that the sum of all angle differences is a multiple of $2 \pi$. Stability of these equilibria can be assessed by examining the Hessian of $V_{\mathrm{circ}}$. This leads to the conclusion that each configuration with $\vert\theta_j-\theta_k\vert = \theta_0 < \tfrac{\pi}{2}$ $\forall (j,k) \in \mathcal{E}$ is locally asymptotically stable under (\ref{eq:3:CircleSynchCT}) or (\ref{eq:3:circlesynchDT}), and all other configurations are unstable. Thus in a stable configuration, consecutive agents in the ring graph are separated by $\theta_0$ on the circle, for some $\theta_0 \in (-\tfrac{\pi}{2},\tfrac{\pi}{2})$ satisfying $N \theta_0 = 2 a \pi$ with $a \in \mathbb{Z}$. The case $\theta_0 = 0$ corresponds to synchronization. When $N \geq 5$, stable configurations exist with $\theta_0 \neq 0$, see Figure \ref{fig:Splayring}. For all these configurations, $\sum_{k=1}^N \; e^{i \theta_k} = 0$, therefore they are said to be \emph{balanced}.

\paragraph*{Stable configurations are graph dependent:} It is currently an open question to characterize, with graph-theoretic properties, which graphs admit no local minima of $V_{\mathrm{circ}}$ different from synchronization. The following result shows that in fact, any configuration that is sufficiently ``spread'' on the circle is stable under the synchronization algorithms for a well-chosen weighted digraph.\\

\begin{proposition}\label{prop:3:IFACall}
\emph{Consider a set of $N$ agents distributed on $S^1$ in a configuration $\lbrace \theta_k \rbrace$ such that for every $k$, there is at least one agent located in $(\theta_k, \theta_k + \pi/2)$ and one located in $(\theta_k - \pi/2, \theta_k)$; such a configuration requires $N \geq 5$. Then there exists a positively weighted and strongly connected $\delta$-digraph making this configuration locally exponentially stable under (\ref{eq:3:CircleSynchCT}) with $\alpha > 0$.}
\end{proposition}\vspace{2mm}

\begin{proof}
Choose nonzero weights $a_{jk}$ only for the $\theta_j \in (\theta_k - \pi/2, \theta_k + \pi/2)$, and such that $r_k := \sum_{j=1}^N a_{jk} \, (x_j - x_k)$ is aligned with $x_k=e^{i \theta_k}$.
\end{proof}

For any of the weight choices that locally stabilize specific configurations, synchronization is also exponentially stable --- but thus only locally. The equilibrium configurations of (\ref{eq:3:circlesynchDT}) and (\ref{eq:3:CircleSynchCT}) different from synchronization are formalized in \cite{MY4} as \emph{consensus} configurations. The same paper also considers the related problem of ``spreading'' agents on the circle, which is formalized with the notions of \emph{anti-consensus} and \emph{balancing} configurations.

\begin{figure}[htb]
\begin{center}\setlength{\unitlength}{1mm}
\begin{picture}(150,50)
\put(65,8.5){\includegraphics[width=50mm]{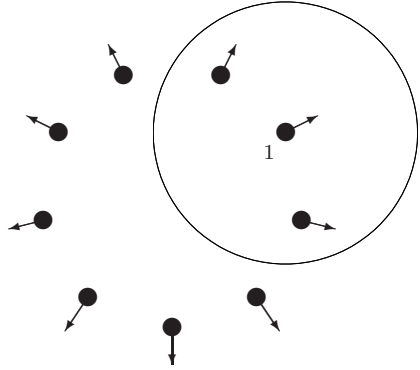}}
\put(88,27.3){{\footnotesize $1$}}
\put(76,5){\circle*{2.5}}
\multiput(64.9,9)(22.2,0){2}{\circle*{2.5}}
\multiput(59,19.2)(34,0){2}{\circle*{2.5}}
\multiput(61.05,30.9)(29.9,0){2}{\circle*{2.5}}
\multiput(69.6,38.45)(12.8,0){2}{\circle*{2.5}}
\put(76,5){\vector(0,-1){5}}
\put(64.9,9){\vector(-2,-3){3}}
\put(87.1,9){\vector(2,-3){3}}
\put(59,19.2){\vector(-4,-1){4.5}}
\put(93,19.2){\vector(4,-1){4.5}}
\put(61.05,30.9){\vector(-2,1){4.2}}
\put(90.95,30.9){\vector(2,1){4.2}}
\put(69.6,38.45){\vector(-1,2){2}}
\put(82.4,38.45){\vector(1,2){2}}
\end{picture}
\end{center}
\caption{Initial conditions for divergent behavior in Vicsek model. The black disks denote agent positions, arrows denote orientations; the circle represents the sensing region of agent $1$.}\label{fig:Vicsek1}
\end{figure}

\paragraph*{Structurally stable divergent behavior in Vicsek model:} The above example of stable equilibria for the undirected ring allows to illustrate a situation where Vicsek model diverges. Consider $N\geq 5$ agents initiated as in Figure \ref{fig:Vicsek1}: (i) initial positions $r_k(0) \in \mathbb{R}^2$ are regularly distributed on a circle such that each agent can sense only its immediate neighbor on the left and on the right; (ii) initial orientations $\theta_k(0) \in S^1$ point radially outwards of the circle formed by the positions. Then the update equation for agent orientations $\theta_k(t)$ is exactly in a stable configuration different from synchronization under a ring interconnection graph. The agents move radially outwards; at a particular time step, all communication links drop. Stability of the equilibrium for the orientations and the fact that all communication links still disconnect at the same instant when positions are slightly shifted ensures that the divergent behavior is observed in an open neighborhood of initial conditions around this ideal situation. Examining stable equilibria of (\ref{eq:3:circlesynchDT}) and (\ref{eq:3:CircleSynchCT}) for $\G$ different from an undirected ring, one sees that the stable divergent behavior remains if the sensing regions are increased such that each agent initially has several neighbors on the left and on the right.

\subsection{Limit sets different from equilibrium}
\label{ss:I:openquestions}

Section \ref{ss:ConvToEq} lists cases where (\ref{eq:3:CircleSynchCT}) or (\ref{eq:3:circlesynchDT}) do not converge to \emph{synchronization}, but still to a set of equilibria. There are also cases where the agents do not converge to a set of equilibria.\\

For fixed undirected graphs, the swarm is ensured to converge to a set of equilibria, except for the discrete-time algorithm when $\beta$ is too small (see Proposition \ref{s5}.\ref{prop:3:DTSynch}). In the latter case, behavior of (\ref{eq:3:circlesynchDT}) is not as clear and the system may run into a limit cycle; see \cite{MY1} or \cite{MYthesis} for details.\\

Periodic and quasi-periodic behaviors can easily be constructed for (\ref{eq:3:CircleSynchCT}) with fixed \emph{directed} graphs.

The simplest such behavior is called \emph{cyclic pursuit}: each agent $k$ is attracted by its neighbors to move (say) clockwise, and the agents keep turning without synchronizing. A classical situation of stable cyclic pursuit is a directed ring graph with consecutive agents separated by $\tfrac{2\pi}{N}$.

In basic cyclic pursuit, agents keep moving on the circle but relative positions remain constant. A more meaningful periodic behavior occurs when \emph{relative positions periodically vary} in time. Such situations can be built with agents partitioned into two sets such that each set is in cyclic pursuit at a different velocity. Start for instance with two unweighted directed ring graphs of $N_1$ and $N_2$ agents, where $N_1 + N_2 = N$ and $N_1 \neq N_2$; consecutive agents in each ring are separated by $\tfrac{2\pi}{N_1}$ and $\tfrac{2\pi}{N_2}$. Then the resulting behavior is satisfactory, but the overall graph is not connected. To obtain a strongly connected graph, each agent of the first ring can be coupled to all the agents of the second ring and conversely; indeed, for a set of regularly spaced agents $\sum_k \; e^{i \theta_k} = 0$, so coupling an agent to such a set does not change its behavior.

Likewise, a \emph{quasi-periodic variation} of relative positions is obtained when several sets of agents move in cyclic pursuit with irrational velocity ratios. This can be built for instance with unitary graph weights and $\alpha=1$, if one set has $x$ agents in a splay state for an undirected ring graph $\Leftrightarrow$ $\tfrac{d}{dt} \theta_k = 0$, a second set has $6$ agents in cyclic pursuit with a directed ring graph $\Leftrightarrow$ $\tfrac{d}{dt} \theta_k = \sqrt{3}$, and a third set has $12$ agents in cyclic pursuit with a directed ring graph $\Leftrightarrow$ $\tfrac{d}{dt}\theta_k = 1$.

\begin{figure}
\begin{center}
\includegraphics[width=120mm]{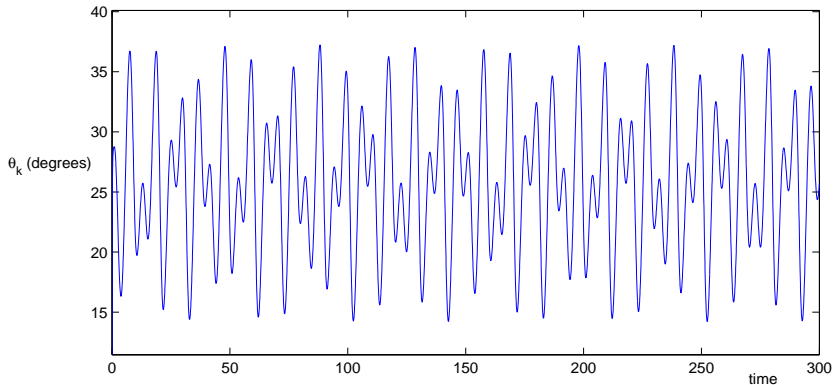}
\end{center}
\caption{Illustration of disorderly-looking quasi-periodic motion of an agent $k$ for which $\tfrac{d}{dt}\theta_k = \sum_{j=1}^3 \; \sin(\theta_j-\theta_k)$ where $\theta_1$ belongs to a regularly spaced undirected ring, $\tfrac{d}{dt}\theta_1 = 0$; $\theta_2$ is in cyclic pursuit with 5 other agents, $\tfrac{d}{dt}\theta_2 = \sqrt{3}$; $\theta_3$ is in cyclic pursuit with 11 other agents, $\tfrac{d}{dt}\theta_3 = 1$.}\label{fig:Chaos1}
\end{figure}

Finally, an example of \emph{disorderly-looking quasi-periodic motion} can be built by adding to the previous situation an agent which is influenced by one agent in each of the three rings; the motion of this agent is illustrated on Figure \ref{fig:Chaos1}.

\begin{figure}
\begin{center}\setlength{\unitlength}{1mm}
\begin{picture}(100,180)
\put(3,0){\includegraphics[width=95mm]{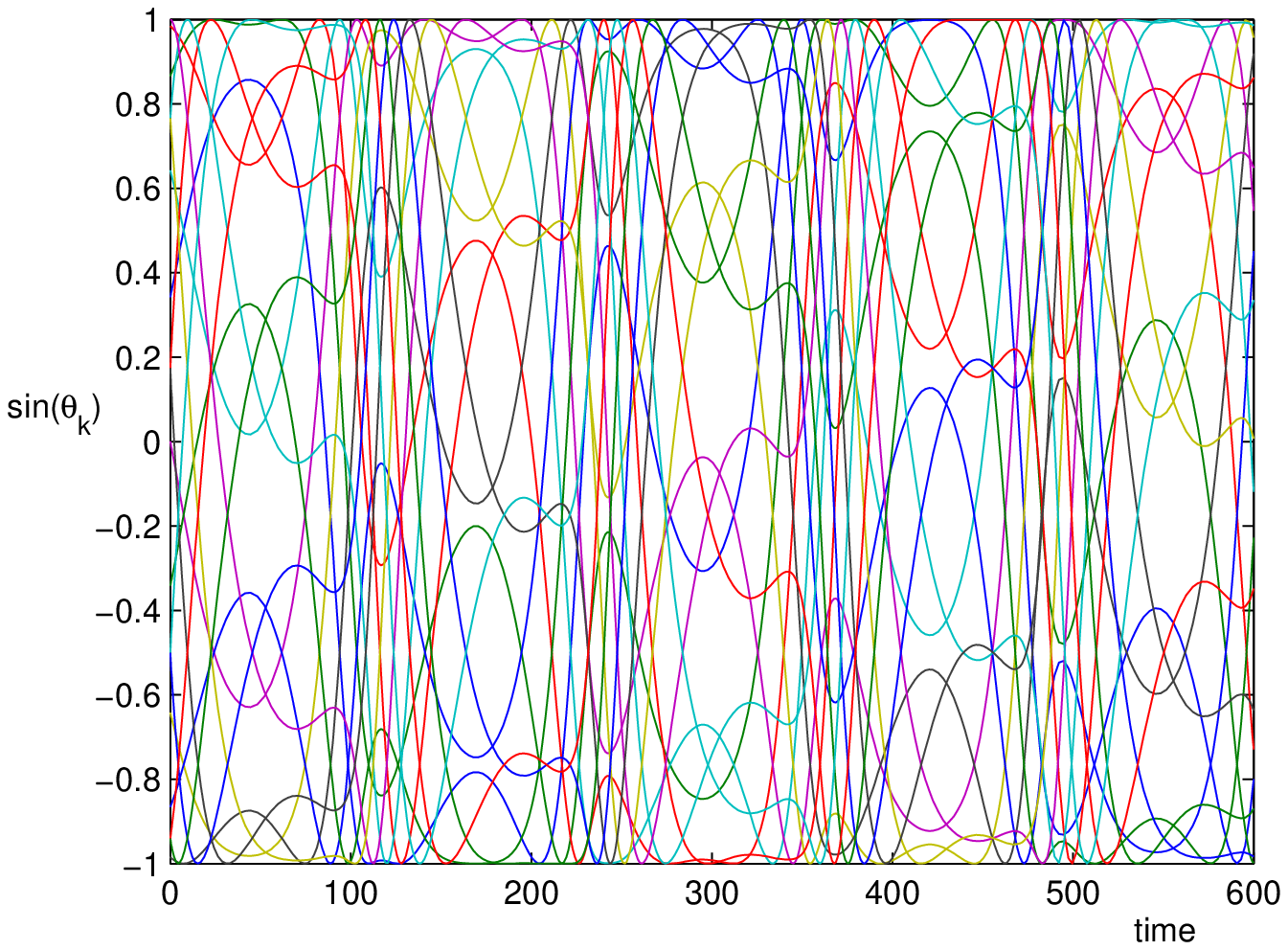}}
\put(0,80){\includegraphics[width=95mm]{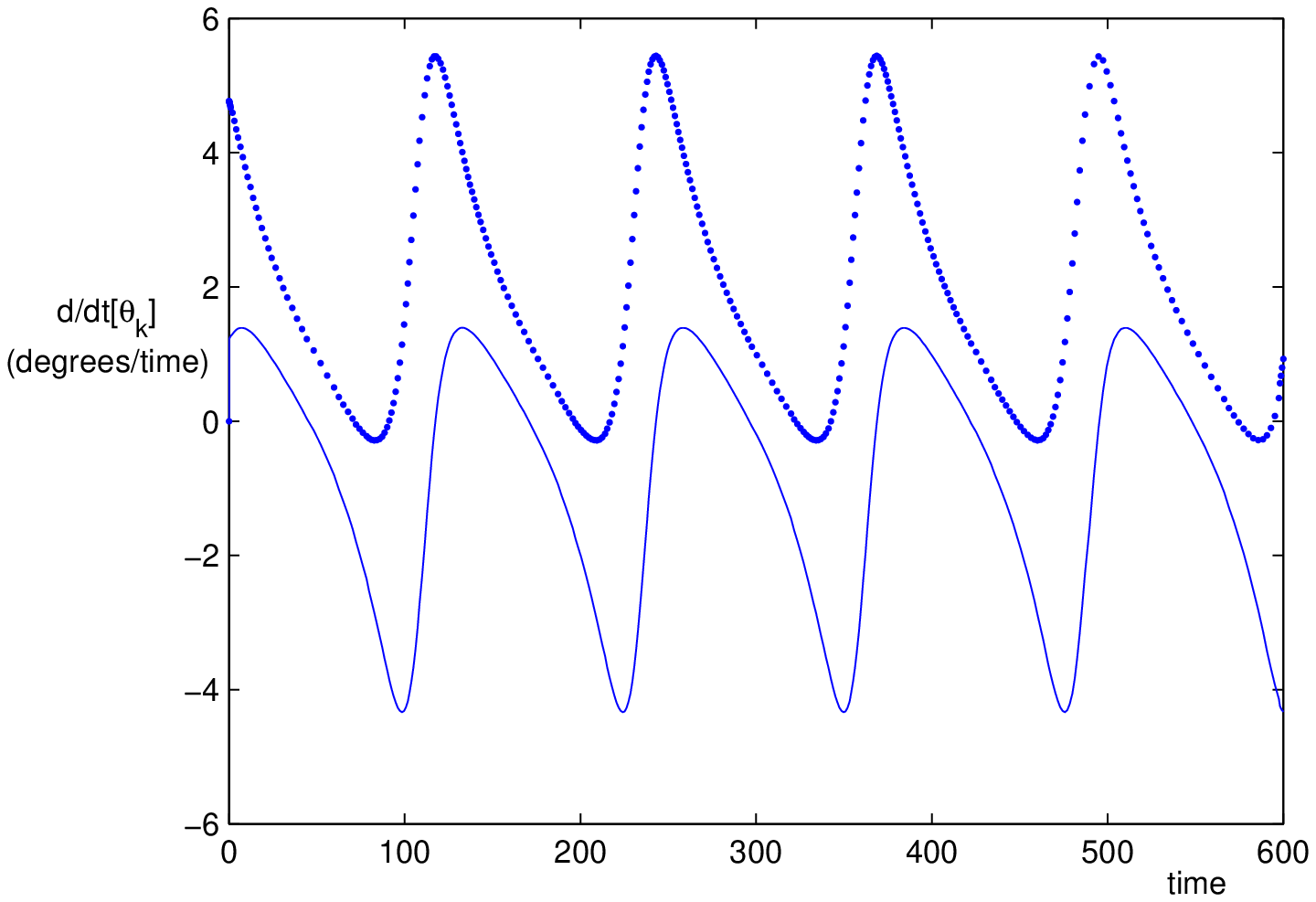} }
\end{picture}
\end{center}
\caption{Motion of agents applying (\ref{eq:3:CircleSynchCT}) with fixed directed coupling, such that they all periodically revert their direction of motion. Set $A$ has $9$ agents regularly spaced by $\tfrac{2\pi}{9}$ on the circle. Set $B$ has $9$ agents also regularly spaced, but initially rotated by $\tfrac{\pi}{18}$ with respect to $A$. In $A$, $\tfrac{d}{dt}\theta_k = 0.04 \, \sin(\theta_j-\theta_k) + 0.05 \sin(\theta_l-\theta_k)$ where $j$ is the agent of $A$ for which $\theta_j-\theta_k = \tfrac{-2\pi}{9}$, and $l$ is the agent of $B$ for which initially $\theta_l-\theta_k = \tfrac{7\pi}{18}$. In $B$, $\tfrac{d}{dt}\theta_k = 0.07 \, \sin(\theta_j-\theta_k) + 0.05 \sin(\theta_l-\theta_k)$ where $j$ is the agent of $B$ for which $\theta_j-\theta_k = \tfrac{2\pi}{9}$ and $l$ is the agent of $A$ for which initially $\theta_l-\theta_k = \tfrac{5\pi}{18}$. Top: velocities of the two sets (continuous curve for $A$, dotted curve for $B$). Bottom: evolution of $\sin(\theta_k)$ for all agents $k$.}\label{fig:Chaos2}
\end{figure}

It is possible to build situations of fixed directed coupling with even more surprising behavior. Figure \ref{fig:Chaos2} represents the motion of two sets of agents in cyclic pursuit, with coupling among agents of the two sets and initial positions such that all the agents periodically revert their direction of motion.\\

\vspace{3mm}

For time-varying graphs, the situation is even more complicated. Since many different configurations can be stable on the circle depending on $\G$, the swarm can be driven towards different equilibria during longer or shorter time spans, implying no particular characterization of the swarm's behavior if $\G(t)$ can be arbitrary. In practice, synchronization is often eventually observed. This is because synchronization is ensured for connected graph sequences as soon as agents all lie in the same semicircle. But other asymptotic behaviors are possible.\\

The diversified, poorly characterized behavior of (\ref{eq:3:circlesynchDT}) or (\ref{eq:3:CircleSynchCT}) with directed and time-varying $\G$ is in strong contrast with the behavior of the consensus algorithm on vector space, which is fully characterized by Proposition \ref{s4}.\ref{prop:3:Moreau}. In addition, Proposition \ref{s4}.\ref{prop:3:Moreau} can be extended to the case where time delays are present along the communication links (see \cite{olfati}), while the behavior of (\ref{eq:3:circlesynchDT}) or (\ref{eq:3:CircleSynchCT}) under time delays is still under investigation even for fixed undirected $\G$. Even for the complete graph, delays may lead to stable synchronized solutions, stable ``spread'' solutions, as well as periodic oscillations \cite{DelayKuramoto}.\\

Finally, it must be noted that the graph modeling inter-agent communication often depends on the states of the agents, like for instance in full Vicsek model \cite{VICSEK}. Studying the interaction of state-dependent graphs with algorithms that are not specifically designed for particular graph behavior is currently a difficult problem, which goes beyond the scope of the present paper.

\section{Algorithms for global synchronization}\label{s6}

Section \ref{s5} highlights that consensus algorithms on the circle may exhibit complicated behaviors. In a control framework, a natural question is then whether the update rules (\ref{eq:3:circlesynchDT}) and (\ref{eq:3:CircleSynchCT}) can be modified to enforce better synchronization properties.

\subsection{Modified coupling for fixed undirected graphs}\label{s:I:newprofile}

For fixed undirected graphs, algorithms (\ref{eq:3:circlesynchDT}) and (\ref{eq:3:CircleSynchCT}) guarantee convergence to an equilibrium set, but stable equilibria different from synchronization may exist. Such spurious equilibria can be rendered unstable by reshaping the way agents are attracted towards their neighbors. A continuous-time setting is chosen for convenience; a similar result could be developed in discrete-time. For simplicity, $\G$ is assumed unweighted.

Consider a continuous-time synchronization algorithm on $S^1$ of the form
\begin{equation}\label{def:profile}
\tfrac{d}{dt} \theta_k = \sum_{j \rsa k} f(\theta_j-\theta_k) \; , \quad k=1,2,...,N \, .
\end{equation}
The function $f : S^1 \rightarrow \mathbb{R}$ is called the \emph{coupling function}. For (\ref{eq:3:CircleSynchCT}), $f(\theta) = \sin(\theta)$. Section \ref{ss:ConvToEq} examines local equilibria of (\ref{eq:3:CircleSynchCT}) for the undirected ring graph and concludes that a configuration is stable if interconnected agents are closer than $\tfrac{\pi}{2}$ (because $\tfrac{d}{d\theta}f(\theta) > 0$ for $\vert \theta \vert < \tfrac{\pi}{2}$), and unstable if they are further apart than $\tfrac{\pi}{2}$ (because $\tfrac{d}{d\theta}f(\theta) < 0$ for $\vert \theta \vert > \tfrac{\pi}{2}$). If $f(\theta)$ is modified to have a positive slope only up to $\tfrac{\pi}{a}$ for some $a > 2$, then connected agents must be closer than $\tfrac{\pi}{a}$ at a stable equilibrium for the ring graph; taking $\tfrac{a}{N} > \tfrac{1}{2}$, it becomes impossible to distribute the agents as on Figure \ref{fig:Splayring} and synchronization is the only stable equilibrium. This motivates the following.

\begin{figure}[htb]
\begin{center}\setlength{\unitlength}{1mm}
\begin{picture}(120,30)
\put(60,0){\vector(0,1){27}}
\put(0,10){\vector(1,0){120}}
\put(60.3,7.7){{\footnotesize{$0$}}}
\multiput(5,0)(50,0){3}{\line(1,2){10}}
\multiput(15,20)(50,0){2}{\line(2,-1){40}}
\put(61,25){$g(\theta)$} \put(118,6){$\theta$}
\multiput(10,9)(25,0){5}{\line(0,1){2}}
\put(9,7){{\footnotesize{$-2\pi$}}}
\put(33,7){{\footnotesize{$-\pi$}}}
\put(84,7){{\footnotesize{$\pi$}}}
\put(110,7){{\footnotesize{$2\pi$}}}
\multiput(0,15)(117,0){2}{$...$}
\multiput(55,9)(10,0){2}{\line(0,1){2}}
\put(65,7){{\footnotesize{$\tfrac{\pi}{N}$}}}
\put(50,7){{\footnotesize{$\tfrac{-\pi}{N}$}}}
\put(54.5,18.5){{\footnotesize{$\tfrac{a \,\pi}{N}$}}}
\put(59,20){\line(1,0){2}}
\end{picture}
\end{center}
\caption{The alternative coupling function $g$.}\label{fig:altprofile}
\end{figure}
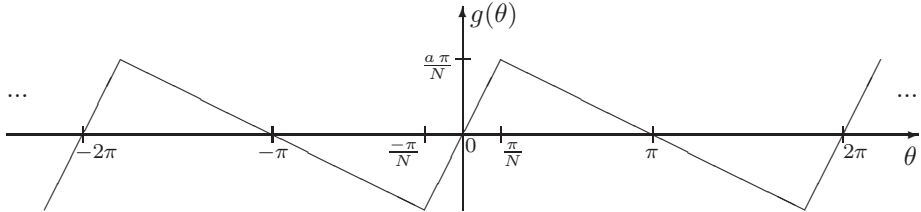

Assume (a bound on) the number $N$ of agents in the swarm is available. Define
\begin{equation}\label{eq:4:profile}
g(\theta) = \left\lbrace \begin{array}{ll}
\frac{-a}{N-1}(\pi+\theta) & \text{ for } \theta \in [-\pi,-\tfrac{\pi}{N}]\\[2mm]
a \, \theta & \text{ for } \theta \in [-\tfrac{\pi}{N},\tfrac{\pi}{N}]\\[2mm]
\frac{a}{N-1}(\pi-\theta) & \text{ for } \theta \in [\tfrac{\pi}{N},\pi]
\end{array}\right.
\end{equation}
for some $a > 0$, extended $2\pi$-periodically outside the above intervals, as represented on Figure \ref{fig:altprofile}. Function $g(\theta)$ is the gradient of (with $2\pi$-periodic extension)
$$(z(\theta))^2 = \left\lbrace \begin{array}{ll}
\tfrac{a \pi^2}{2 N (N-1)} + \tfrac{a}{N-1}(-\pi\, \theta - \tfrac{\theta^2}{2}) & \text{ for } \theta \in [-\pi,\, -\tfrac{\pi}{N}]\\[2mm]
\tfrac{a}{2} \theta^2 & \text{ for } \theta \in [-\tfrac{\pi}{N},\, \tfrac{\pi}{N}]\\[2mm]
\tfrac{a \pi^2}{2 N (N-1)} + \tfrac{a}{N-1}(\pi\, \theta - \tfrac{\theta^2}{2}) & \text{ for } \theta \in [\tfrac{\pi}{N},\, \pi]
\end{array}
\right.$$
which is even, has a minimum for $\theta=0$, a maximum for $\theta = \pi$ and evolves monotonically and continuously in between, similarly to the sinusoidal distance measure $\left( 2 \sin(\tfrac{\theta}{2}) \right)^2$.

Choosing $f(\theta)=g(\theta)$ in (\ref{def:profile}) defines a synchronization algorithm which satifies all invariance and communication constraints and whose only stable configuration is synchronization for fixed undirected $\G$. Rigorously, the edges in $g(\theta)$ should be smoothed to make it continuously differentiable everywhere; this changes nothing to the general argument.\\

\begin{proposition}\label{prop:4:alt.profile}
\emph{Consider a swarm of $N$ agents, interconnected according to a connected fixed undirected graph $\G$, that evolve on $S^1$ by applying (\ref{def:profile}) with $f(\theta)=g(\theta)$ defined by (a smoothed version of) (\ref{eq:4:profile}). The agents always converge to the set of equilibria corresponding to the critical points of $V_g = \tfrac{1}{2}\, \sum_k \sum_{j\rsa k}\; \left( z(\theta_j - \theta_k) \right)^2$. Moreover, the only asymptotically stable equilibrium is synchronization.}
\end{proposition}\vspace{2mm}

\begin{proof}
(see \cite{MYthesis} for a full proof) The agents always converge to a set of critical points of $V_g$ because the algorithm is a gradient descent for $V_g$. Synchronization, as the global minimum of $V_g$, is stable. Stability of other equilibria is characterized by examining the Hessian of $V_g$. It is shown that for the interaction function $g(\theta)$, if the Hessian is positive semidefinite with $0$ eigenvalue only in the direction of uniform motion ($\theta_k \rightarrow \theta_k + a$ $\forall k \in \mathcal{V}$), then the graph $\G_p$, containing edge $\lbrace j ,k \rbrace$ if and only if $j$ and $k$ are closer than $\tfrac{\pi}{N}$, must be connected. Then all agents must be within a semicircle and convexity arguments impose synchronization like on the real line.
\end{proof}

\subsection{Introducing randomness in link selection}
\label{s:I:randomlinks}

The modified coupling function solves the problem of spurious local minima in $V_{\mathrm{circ}}$. However, it may introduce numerous unstable equilibria. Moreover, for varying and directed graphs, the behavior of (\ref{def:profile}) is not better characterized than for (\ref{eq:3:CircleSynchCT}). The present section introduces a so-called ``Gossip Algorithm'' (see \cite{BoydGossip} and references therein) in order to improve synchronization behavior on the circle. Thanks to the introduction of randomness, it achieves \emph{global asymptotic synchronization with probability $1$} for directed and time-varying $\G$. It is described in discrete-time for easier formulation.

The nice convergence properties of (\ref{eq:3:circlesynchDT}) and (\ref{eq:3:CircleSynchCT}) when $\G$ is a tree motivate to keep the update law (\ref{eq:3:circlesynchDT}), but at each time select at most one of the in-neighbors in $\G(t)$ for each agent. In order to satisfy the equivalence of all agents, $k$ may not privilege any of its neighbors --- it is just allowed, for weighted $\G$, to take the different weights of the corresponding edges into account. Always choosing the neighbor with maximum weight could disconnect the swarm. Therefore, it is necessary to select the retained neighbor \emph{randomly} among the $j \rsa k$. A natural probability distribution for neighbor selection would follow the weights of the edges, but in theory any distribution with nonzero weight on each link is admissible. With the proposed edge selection procedure, the neighbor chosen at time $t+1$ is independent of the neighbor chosen at time $t$ (up to, for varying graphs, a possible dependence of $\G(t+1)$ on $\G(t)$).\vspace{2mm}

\noindent \textbf{Gossip algorithm (directed).} At each update $t$,
\begin{enumerate}[topsep=0mm, parsep=0mm, itemsep=0mm]
\item[1.] each agent $k$ randomly selects an agent $j \rightsquigarrow k$ with probability $a_{jk} / \left( \beta + \sum_{l \rsa k} \, a_{lk} \right)$, where $\beta > 0$ is the weight for choosing no agent;
\item[2.] $\theta_k(t+1) = \theta_j(t)$ if agent $k$ chooses neighbor $j$ at time $t$, and $\theta_k(t+1) = \theta_k(t)$ if it chooses no neighbor.
\end{enumerate}\vspace{2mm}

A variant of the Gossip Algorithm exists in which the random graph is \emph{undirected} at each time step. In this variant, agents have to select each other in order to move by averaging their positions. The advantages of this variant on $S^1$ are not clear, and proving global asymptotic synchronization with probability $1$ is somewhat more difficult; see \cite{MY9} for details.

The directed Gossip Algorithm proposed above is extreme in the sense that agent $k$ directly jumps to the position of its selected neighbor. A more moderate directed Gossip Algorithm would apply the update law $\theta_k(t+1)=\mathrm{arg}(\alpha \, e^{i\theta_k(t)} + e^{i\theta_j(t)})$ with $\alpha > 0$.\\

The authors of \cite{BoydGossip} perform a detailed analysis of a Gossip Algorithm for synchronization in vector spaces. Convergence towards synchronization is always ensured on vector spaces and the problem is to quantify the convergence \emph{rate} as a function of $\G$ and probability (i.e. weights) distribution. On the circle, convergence towards synchronization is not obvious a priori, but in fact it holds under the same assumptions.\\

\begin{definition}\label{def:4:probability.synch}
\emph{$N$ agents \emph{asymptotically synchronize with probability $1$} if for any initial condition, for any $\varepsilon > 0$ and $\kappa \in (0,1)$, there exists a time $T$ after which the maximal distance $\vert \theta_k(T) - \theta_j(T) \vert$ between any pair of agents is smaller than $\varepsilon$ with probability larger than $\kappa$.}
\end{definition}\\

\begin{proposition}\label{prop:4:GossipDir}
\emph{Consider a set of $N$ agents interconnected according to a uniformly connected $\delta$-digraph $\G$. If the agents apply the (directed) Gossip Algorithm with a fixed finite $\beta > 0$, then they asymptotically synchronize with probability $1$.}
\end{proposition}\vspace{2mm}

\begin{proof}
(see \cite{MY9} for a full proof) Consider a set of links forming a directed tree of root $k$. For $\G(t)$ uniformly connected, there exists a sequence of link choices, implementable on any time interval $[t, t+T_s]$ for some finite $T_s$, that builds this tree (and only this tree, potentially selecting no move for many time instants!) sequentially from the root to its leaves. This sequence synchronizes the agents at the position of the root \emph{for any initial conditions}. Moreover, its finite length implies a finite (though potentially tiny) probability to be selected during any time interval of length $T_s$. Therefore when $t$ goes to $+\infty$, the probability that this sequence never appears goes to $0$. Since it suffices that the sequence appears once to ensure synchronization, this concludes the proof.
\end{proof}

The convergence proof of Proposition \ref{s6}.\ref{prop:4:GossipDir} can be adapted for the moderate version of the directed Gossip Algorithm, with inertia $\alpha > 0$, as described in \cite{MYthesis}.

Under the initial Gossip Algorithm (inertia-less, i.e. with $\alpha=0$), agents in fact jump between a discrete set of possible positions corresponding to the initial positions of the $N$ agents. This highlights that the directed Gossip Algorithm can in fact be applied on \emph{any set of symbols}. Proposition \ref{s6}.\ref{prop:4:GossipDir} purely relies on the evolution of agents between $N$ different ``symbols'', \emph{completely independently of the underlying manifold}. Every time a position is left empty (implying that the synchronization process progresses as an agent joins other ones), that position can never be reached again in the future; this process goes on until all agents are on the same position after a finite time.

In this context, a natural measure of convergence rate is the \emph{expected synchronization time}, i.e. the average time, over all possible link choices, after which all agents are on the same position. A Markov chain framework can be applied to obtain an explicit formula for the expected synchronization time, at least for fixed graph $\G$, see \cite{MYthesis}; unfortunately, the complexity of the explicit formula grows exponentially with the number of agents, and it seems difficult to extract the influence of the link choice probability distribution. The expected synchronization time --- including its numerical value --- is independent of $S^1$ and \emph{independent of the initial positions of the agents} (unless some agents are initially perfectly synchronized). The only remaining parameters are the graph $\G$ and the probability distribution. There seems to be an interesting interplay between these two parameters; it is therefore challenging for the agents to optimize their convergence rate based on their local information.

Simulations confirm that the Gossip Algorithm favors global synchronization. However, they also highlight that the probabilistic setting can lead to unnecessarily slow convergence rates: often a set of agents partitioned into two groups, located at $\theta_a$ and $\theta_b$ respectively, keep oscillating between these two positions for an appreciable time before they synchronize on one of them.

\subsection{Algorithms using auxiliary variables}
\label{s:I:aux}

In the synchronization framework, the non-convexity of $S^1$ can be ``cheated'' if the agents are able to communicate \emph{auxiliary variables} in addition to their positions on the circle. A different viewpoint on this procedure is that the limited number of communication links for information flow is compensated by sending larger communication packets along existing links. Such strategies allow to recover the synchronization properties of vector spaces for almost all initial conditions. Their potential interest lies more in engineering applications than in physical modeling, where communication of auxiliary variables is questionable. In the present paper, the relevant algorithms are just briefly mentioned for completeness.\\

A first use of auxiliary variables is to build a reduced communication network, in which a leader is identified which then attracts all the agents. The leader election / spanning tree construction would be completely independent of the agents' motion on $S^1$ and should be achieved after a finite time; then applying a synchronization algorithm with this leader / spanning tree would ensure synchronization. See \cite{SpanningTree} for a distributed algorithm that achieves this preliminary network construction. Note that for termination of the network construction after finite time, (a bound on) the total number $N$ of agents must be available to each agent. Also, it is not clear how well the spanning tree construction can be adapted to time-varying graphs.\\

A second use of auxiliary variables is to reach consensus on a \emph{reference synchronization point}: an auxiliary variable in $\mathbb{R}^2$ is associated to each agent, and a consensus algorithm is run on these auxiliary variables. Each agent individually tracks the projection of its auxiliary variable on $S^1$. Thanks to the fact that $S^1$ is equivalent to the Lie group $SO(2)$, this can be implemented in a way that satisfies the coupling symmetry hypotheses. For uniformly connected $\G$, the consensus algorithm on auxiliary variables defines a point in the plane, and synchronization is ensured for almost all initial conditions. See \cite{MY2} for details and a formal proof.

\section{Generalizations on compact homogeneous manifolds}\label{s7}

Although the discussion in this paper focuses on the circle for simplicity, it is representative of more general manifolds. Compact homogeneous manifolds are manifolds that can be viewed as the quotient of a compact Lie group by one of its subgroups. They include the $n$-dimensional sphere $S^n$, all compact Lie groups like e.g. the group of rotations $SO(n)$, and the Grassmann manifolds $\Gr$.

In \cite{MY4}, the developments of Section \ref{s3} are extended to compact homogeneous manifolds $\mathcal{M}$ that are embedded in $\mathbb{R}^m$ such that the Euclidean norm $\Vert x \Vert$ in $\mathbb{R}^m$ is constant over $x \in \mathcal{M}$. The cost function and gradient algorithms are generalized, and an intepretation in terms of moving towards an appropriately defined \emph{average} of positions on $\mathcal{M}$ is given. In particular, the cost function measures the Euclidean distance between agents \emph{in the embedding space}, also called the \emph{chordal distance}. Convergence properties are analyzed, and stable local equilibria are formalized as \emph{consensus} configurations. An ``opposite'' algorithm, in which agents move away from their neighbors to spread on the circle, is proposed along the same lines and its convergence properties are formalized with \emph{anti-consensus} and \emph{balancing} configurations.

Regarding global synchronization, the results of Section \ref{s4} for trees and complete graphs as well as agents initially located within convex sets remain valid on more general manifolds. Obstacles like those illustrated in Section \ref{s5} appear as well. The modified coupling algorithm of Section \ref{s:I:newprofile} should, in principle, have an extension on general manifolds, but the current convergence proof is algebraic rather than geometric and cannot be repeated. The (directed) Gossip Algorithm of Section \ref{s:I:randomlinks} works on any manifold --- in fact even any set; the undirected variant is less easily generalized. The algorithm mentioned in Section \ref{s:I:aux}, specifying a ``reference synchronization state'' with auxiliary variables, is generalized to compact homogeneous manifolds in \cite{MY4} (modulo the fact that a meaningful communication of auxiliary variables sometimes requires a common external reference frame). The same type of algorithm with auxiliary variables can also be used for balancing on the circle and on compact homogeneous manifolds.

\section{Conclusion}

The present paper proposes a natural extension of the ``consensus'' framework, where agents have to agree on a common value, from values in \emph{vector spaces} to values on the \emph{circle}. This leads to a related interpretation of Kuramoto and Vicsek models of the literature.

Although convergence is similar to vector spaces in some specific situations, in general the behavior of consensus algorithms on the circle is much more diversified, allowing local equilibria, limit cycles, and essentially any type of behavior when the interconnection graph can vary freely.

Global synchronization properties can be recovered on the circle (i) by modifying the coupling between agents (Section \ref{s:I:newprofile}), (ii) in a stochastic setting where at most one neighbor is randomly selected at each time step (Section \ref{s:I:randomlinks}), and (iii) less surprisingly, by assisting the agreement process with auxiliary variables in a vector space (Section \ref{s:I:aux}).

Several questions remain open to characterize and optimize the \emph{global} behavior of these simple consensus algorithms on the circle; this indicates possible directions of interest for future research.


\section*{Acknowledgments}

This paper presents research results of the Belgian Network DYSCO (Dynamical Systems, Control, and Optimization), funded by the Interuniversity Attraction Poles Programme, initiated by the Belgian State, Science Policy Office. The scientific responsibility rests with its authors. Drs. L. Scardovi and J. Hendrickx as well as Profs. V. Blondel, E. Tuna, P-A. Absil and N. Leonard are acknowledged for interesting discussions related to this subject. The second author is supported as an FNRS fellow (Belgian Fund for Scientific Research).

\section*{Appendix: Notions of graph theory}

In the framework of coordination with limited interconnections between agents, it is customary to represent communication links by means of a \emph{graph} (see for instance \cite{Graphs2,Graphs}).\\

\begin{definition} \label{def:1:digraph} \emph{A \emph{directed graph} $\G(\mathcal{V},\mathcal{E})$ (short \emph{digraph} $\G$) is composed of a finite set $\mathcal{V}$ of \emph{vertices}, and a set $\mathcal{E}$ of \emph{edges} which represent interconnections among the vertices as ordered pairs $(j,k)$ with $j$ and $k \in \mathcal{V}$.}

\emph{A \emph{weighted digraph} $\G(\mathcal{V},\mathcal{E},\mathcal{A})$ is a digraph associated with a set $\mathcal{A}$ that assigns a positive \emph{weight} $a_{jk} \in \mathbb{R}_{>0}$ to each edge $(j,k) \in \mathcal{E}$.}\\
\end{definition}

An unweighted graph is often considered as a weighted graph with unit weights. A weighted graph can be defined by its vertices and weights only, by extending the weight set to all pairs of vertices and imposing $a_{jk} = 0$ if and only if $(j,k)$ does not belong to the edges of $\G$. A digraph is said to be \emph{undirected} if $a_{jk} = a_{kj}$ $\forall j,k \in \mathcal{V}$. It may happen that $(j,k) \in \mathcal{E}$ whenever $(k,j) \in \mathcal{E}$ $\forall j,k \in \mathcal{V}$, but $a_{jk} \neq a_{kj}$ for some $j,k \in \mathcal{V}$; in this case the graph is called \emph{bidirectional}. Equivalently, an unweighted undirected graph can be defined as a set of vertices and a set of \emph{unordered} pairs of vertices.

In the present paper, each agent is identified with a vertex of a graph; the $N$ agents = vertices are designed by positive integers $1,2,...,N$, so $\mathcal{V} = \lbrace 1,2,...,N \rbrace$. The presence of edge $(j,k)$ has the meaning that agent $j$ sends information to agent $k$, or equivalently, agent $k$ measures quantities concerning agent $j$. It is assumed that no ``communication link'' is needed for an agent $k$ to get information about itself, so $\G$ contains no self-loops: $(k,k) \notin \mathcal{E}$ $\forall k \in \mathcal{V}$. In visual representations of a graph, a vertex is depicted by a point, and edge $(j,k)$ by an arrow from $j$ to $k$. Therefore a frequent alternative notation for $(j,k) \in \mathcal{E}$ is $j \rsa k$. One also says that $j$ is an \emph{in-neighbor} of $k$ and $k$ is an \emph{out-neighbor} of $j$. In the visual representation of an undirected graph, all arrows are bidirectional; therefore arrowheads are usually dropped. One simply says that $j$ and $k$ are \emph{neighbors} and writes $j \sim k$ instead of $j \rsa k$ and $k \rsa j$. The \emph{in-degree} of vertex $k$ is $d^{(i)}_k = \sum_{j=1}^N a_{jk}$. The \emph{out-degree} of vertex $k$ is $d^{(o)}_k = \sum_{j=1}^N a_{kj}$. A digraph is said to be \emph{balanced} if $d^{(i)}_k = d^{(o)}_k$ $\forall k \in \mathcal{V}$; in particular, undirected graphs are balanced.\\

The \emph{adjacency matrix} $A \in \mathbb{R}^{N \times N}$ of a graph $\G$ contains $a_{jk}$ in row $j$, column $k$; it is symmetric if and only if $\G$ is undirected. Denote by $\vert \mathcal{E} \vert$ the number of edges in $\G$. For a digraph $\G$, each column of the \emph{incidence matrix} $B \in (\lbrace -1, 0, 1 \rbrace)^{N \times \vert \mathcal{E} \vert}$ corresponds to one edge and each row to one vertex; if column $m$ corresponds to edge $(j,k)$, then
$$b_{jm} = -1 \, , \qquad b_{km} = 1 \qquad \text{and} \qquad b_{lm} = 0 \text{ for } l \notin \lbrace j,k \rbrace \, . $$
For an undirected graph $\G$, each column corresponds to an undirected edge; an \emph{arbitrary orientation} $(j,k)$ or $(k,j)$ is chosen for each edge and $B$ is built for the resulting directed graph. Thus $B$ is not unique for a given $\G$, but $\G$ is unique for a given $B$.

The in- and out-degrees of vertices $1,2,...,N$ can be assembled in diagonal matrices $D^{(o)}$ and $D^{(i)}$. The \emph{in-Laplacian} of $\G$ is $L^{(i)} = D^{(i)} - A$. Similarly, the associated \emph{out-Laplacian} is $L^{(o)} = D^{(o)} - A$. For a balanced graph $\G$, the \emph{Laplacian} $L = L^{(i)} = L^{(o)}$. The standard definition of Laplacian $L$ is for undirected graphs. For the latter, $L$ is symmetric and, remarkably, $L=B B^T$. For general digraphs, by construction, $(\mathbf{1}_N)^T \, L^{(i)} = 0$ and $L^{(o)} \, \mathbf{1}_N = 0$ where $\mathbf{1}_N$ is the column vector of $N$ ones. The spectrum of the Laplacian reflects several interesting properties of the associated graph, specially in the case of undirected graphs, see for example \cite{Graphs}. In particular, it reflects its \emph{connectivity} properties.\\

A \emph{directed path} of length $l$ from vertex $j$ to vertex $k$ is a sequence of vertices $v_0, v_1,...,v_l$ with $v_0 = j$ and $v_l = k$ and such that $(v_m, v_{m+1}) \in \mathcal{E}$ for $m=0,1,...,l-1$. An \emph{undirected path} between vertices $j$ and $k$ is a sequence of vertices $v_0, v_1,...,v_l$ with $v_0 = j$ and $v_l = k$ and such that $(v_m, v_{m+1}) \in \mathcal{E}$ or $(v_{m+1}, v_m) \in \mathcal{E}$, for $m=0,1,...,l-1$. A digraph $\G$ is \emph{strongly connected} if it contains a directed path from every vertex to every other vertex (and thus also back to itself). A digraph $\G$ is \emph{root-connected} if it contains a node $k$, called the \emph{root}, from which there is a path to every other vertex (but not necessarily back to itself). A digraph $\G$ is \emph{weakly connected} if it contains an undirected path between any two of its vertices. For an undirected graph $\G$, all these notions become equivalent and are simply summarized by the term \emph{connected}. For $\G$ representing interconnections in a network of agents, clearly coordination can only take place if $\G$ is connected. If this is not the case, coordination will only be achievable separately in each connected component of $\G$. A more interesting discussion of connectivity arises when the graph $\G$ can vary with time. Before discussing this case, the following summarizes some spectral properties of the Laplacian that are linked to the connectivity of the associated graph.\\

\noindent \textbf{Properties (Laplacian):} Consider the out-Laplacian $L^{(o)}$ of digraph $\G$.
\begin{itemize}[topsep=0mm, parsep=0mm, itemsep=0mm]
\item[(a)] All eigenvalues of $L^{(o)}$ have nonnegative real parts.
\item[(b)] If $\G$ is strongly connected, then $\,0$ is a simple eigenvalue of $L^{(o)}$.
\item[(c)] Expression $x^T L x$, with $x \in \mathbb{R}^N$, is positive semidefinite if and only if $\G$ is balanced.
\end{itemize}
If $\G$ is undirected, the Laplacian $L$ has the following properties.
\begin{itemize}[topsep=0mm, parsep=0mm, itemsep=0mm]
\item[(d)] $L$ is symmetric positive semidefinite.
\item[(e)] The algebraic and geometric multiplicity of $\,0$ as an eigenvalue of $L$ is equal to the number of connected components in $\G$.\\
\end{itemize}

In a coordination problem, interconnections among agents can vary with time, as some links are dropped and others are established. In this case, the communication links are represented by a \emph{time-varying graph} $\G (t)$ in which the vertex set $\mathcal{V}$ is fixed (by convention), but edges $\mathcal{E}$ and weights $\mathcal{A}$ can depend on time. All the previous definitions carry over to time-varying graphs; simply, each quantity depends on time. To prevent edges from vanishing or growing indefinitely, the present paper considers \emph{$\delta$-digraphs}, for which the elements of $A(t)$ are bounded and satisfy the threshold $a_{jk}(t) \geq \delta > 0$ $\forall (j,k) \in \mathcal{E}(t)$, for all $t$.
In addition, in continuous-time $\G$ is assumed to be piecewise continuous. For $\delta$-digraphs $\G (t)$, it is intuitively clear that coordination may be achieved if information exchange is ``sufficiently frequent'', without requiring it to take place all the time. The following definition of ``integrated connectivity over time'' can be found in \cite{hendrickx1,MOREAU2,MOREAU,TsitsiklisThesis}.\\

\begin{definition}\label{def:1:uniform.connectivity} (from \cite{MOREAU2,MOREAU}) \emph{In discrete-time, for a  $\delta$-digraph $\G(\mathcal{V},\mathcal{E}(t),\mathcal{A}(t))$ and some constant $T \in \mathbb{Z}_{\geq 0}$, define the graph $\bar{\G}(\mathcal{V},\bar{\mathcal{E}}(t),\bar{\mathcal{A}}(t))$ where $\bar{\mathcal{E}}(t)$ contains all edges that appear in $\G (\tau)$ for $\tau \in [t, t+T]$ and $\bar{a}_{jk}(t) = \sum_{\tau=t}^{t+T} a_{jk}(\tau)$. Similarly, in continuous-time, for a  $\delta$-digraph $\G(\mathcal{V},\mathcal{E}(t),\mathcal{A}(t))$ and some constant $T \in \mathbb{R}_{> 0}$, define the graph $\bar{\G}(\mathcal{V},\bar{\mathcal{E}}(t),\bar{\mathcal{A}}(t))$ by}
\begin{eqnarray*}
& & \bar{a}_{jk}(t) = \left\lbrace 
\begin{array}{ll}
\int_{t}^{t+T} a_{jk}(\tau) d\tau & \text{if } \int_{t}^{t+T} a_{jk}(\tau) d\tau \geq \delta \\
0 & \text{if } \int_{t}^{t+T} a_{jk}(\tau) d\tau < \delta
\end{array} \right. \\
& & (j,k) \in \bar{\mathcal{E}}(t) \text{ if and only if } \bar{a}_{jk}(t) \neq 0 \; .
\end{eqnarray*}
\emph{Then $\G(t)$ is said to be \emph{uniformly connected over $T$} if there exists a time horizon $T$ and a vertex $k \in \mathcal{V}$ such that $\bar{\G}(t)$ is root-connected with root $k$ for all $t$.}\\
\end{definition}

The following graphs are regularly used in the present dissertation.
\begin{itemize}[topsep=0mm, parsep=0mm, itemsep=0mm]
\item The (equally weighted) \emph{complete graph} is an unweighted, undirected graph that contains an edge between any pair of vertices.
\item An \emph{undirected ring} or \emph{cycle graph} on $N>1$ vertices is equivalent to an undirected path containing all vertices, to which is added an edge between the extreme vertices of the path. Similarly, a \emph{directed ring} or \emph{cycle graph} on $N>1$ vertices is equivalent to a directed path containing all vertices, to which is added an edge from the last to the first vertex in the path.
\item An \emph{undirected tree} is a connected undirected graph in which it is impossible to select a subset of at least $3$ vertices and a subset of edges among them to form an undirected cycle. A \emph{directed tree} of root $k$ is a root-connected digraph of root $k$, in which every vertex can be reached from $k$ by following one and only one directed path.

In a directed tree $\G$, the (unique) in-neighbor of a vertex $j$ is called its \emph{parent} and its out-neighbors are its \emph{children}. The root has no parent, and the vertices with no children are called the \emph{leaves}. This can be carried over to an undirected graph after selecting an arbitrary root.\\
\end{itemize}

\bibliographystyle{plain}
\bibliography{asrsBib.bib}

\end{document}